
\magnification=1200

\hsize=14cm
\vsize=20.5cm
\parindent=0cm   \parskip=0pt
\pageno=1

\def\ind{\hskip 1cm\relax}

\hoffset=-0.5cm
\voffset=-1cm

\pretolerance=500 \tolerance=1000  \brokenpenalty=5000

\catcode`\@=11

\font\eightrm=cmr8         \font\eighti=cmmi8
\font\eightsy=cmsy8        \font\eightbf=cmbx8
\font\eighttt=cmtt8        \font\eightit=cmti8
\font\eightsl=cmsl8        \font\sixrm=cmr6
\font\sixi=cmmi6           \font\sixsy=cmsy6
\font\sixbf=cmbx6


\font\tengoth=eufm10       \font\tenbboard=msbm10
\font\eightgoth=eufm8      \font\eightbboard=msbm8
\font\sevengoth=eufm7      \font\sevenbboard=msbm7
\font\sixgoth=eufm6        \font\fivegoth=eufm5

\font\tencyr=wncyr10       
\font\eightcyr=wncyr8      
\font\sevencyr=wncyr7      
\font\sixcyr=wncyr6

\newfam\gothfam           \newfam\bboardfam
\newfam\cyrfam

\def\tenpoint{%
  \textfont0=\tenrm \scriptfont0=\sevenrm
\scriptscriptfont0=\fiverm
  \def\rm{\fam\z@\tenrm}%
  \textfont1=\teni  \scriptfont1=\seveni  \scriptscriptfont1=\fivei
  \def\oldstyle{\fam\@ne\teni}\let\old=\oldstyle
  \textfont2=\tensy \scriptfont2=\sevensy
\scriptscriptfont2=\fivesy
  \textfont\gothfam=\tengoth \scriptfont\gothfam=\sevengoth
  \scriptscriptfont\gothfam=\fivegoth
  \def\goth{\fam\gothfam\tengoth}%
  \textfont\bboardfam=\tenbboard \scriptfont\bboardfam=\sevenbboard
  \scriptscriptfont\bboardfam=\sevenbboard
  \def\bb{\fam\bboardfam\tenbboard}%
 \textfont\cyrfam=\tencyr \scriptfont\cyrfam=\sevencyr
  \scriptscriptfont\cyrfam=\sixcyr
  \def\cyr{\fam\cyrfam\tencyr}%
  \textfont\itfam=\tenit
  \def\it{\fam\itfam\tenit}%
  \textfont\slfam=\tensl
  \def\sl{\fam\slfam\tensl}%
  \textfont\bffam=\tenbf \scriptfont\bffam=\sevenbf
  \scriptscriptfont\bffam=\fivebf
  \def\bf{\fam\bffam\tenbf}%
  \textfont\ttfam=\tentt
  \def\tt{\fam\ttfam\tentt}%
  \abovedisplayskip=12pt plus 3pt minus 9pt
  \belowdisplayskip=\abovedisplayskip
  \abovedisplayshortskip=0pt plus 3pt
  \belowdisplayshortskip=4pt plus 3pt
  \smallskipamount=3pt plus 1pt minus 1pt
  \medskipamount=6pt plus 2pt minus 2pt
  \bigskipamount=12pt plus 4pt minus 4pt
  \normalbaselineskip=12pt
  \setbox\strutbox=\hbox{\vrule height8.5pt depth3.5pt width0pt}%
  \let\bigf@nt=\tenrm       \let\smallf@nt=\sevenrm
  \normalbaselines\rm}

\def\eightpoint{%
  \textfont0=\eightrm \scriptfont0=\sixrm
\scriptscriptfont0=\fiverm
  \def\rm{\fam\z@\eightrm}%
  \textfont1=\eighti  \scriptfont1=\sixi  \scriptscriptfont1=\fivei
  \def\oldstyle{\fam\@ne\eighti}\let\old=\oldstyle
  \textfont2=\eightsy \scriptfont2=\sixsy
\scriptscriptfont2=\fivesy
  \textfont\gothfam=\eightgoth \scriptfont\gothfam=\sixgoth
  \scriptscriptfont\gothfam=\fivegoth
  \def\goth{\fam\gothfam\eightgoth}%
  \textfont\cyrfam=\eightcyr \scriptfont\cyrfam=\sixcyr
  \scriptscriptfont\cyrfam=\sixcyr
  \def\cyr{\fam\cyrfam\eightcyr}%
  \textfont\bboardfam=\eightbboard
\scriptfont\bboardfam=\sevenbboard
  \scriptscriptfont\bboardfam=\sevenbboard
  \def\bb{\fam\bboardfam}%
  \textfont\itfam=\eightit
  \def\it{\fam\itfam\eightit}%
  \textfont\slfam=\eightsl
  \def\sl{\fam\slfam\eightsl}%
  \textfont\bffam=\eightbf \scriptfont\bffam=\sixbf
  \scriptscriptfont\bffam=\fivebf
  \def\bf{\fam\bffam\eightbf}%
  \textfont\ttfam=\eighttt
  \def\tt{\fam\ttfam\eighttt}%
  \abovedisplayskip=9pt plus 3pt minus 9pt
  \belowdisplayskip=\abovedisplayskip
  \abovedisplayshortskip=0pt plus 3pt
  \belowdisplayshortskip=3pt plus 3pt
  \smallskipamount=2pt plus 1pt minus 1pt
  \medskipamount=4pt plus 2pt minus 1pt
  \bigskipamount=9pt plus 3pt minus 3pt
  \normalbaselineskip=9pt
  \setbox\strutbox=\hbox{\vrule height7pt depth2pt width0pt}%
  \let\bigf@nt=\eightrm     \let\smallf@nt=\sixrm
  \normalbaselines\rm}

\tenpoint

\def\pc#1{\bigf@nt#1\smallf@nt}         \def\pd#1 {{\pc#1} }

\def\^#1{\if#1i{\accent"5E\i}\else{\accent"5E #1}\fi}
\def\"#1{\if#1i{\accent"7F\i}\else{\accent"7F #1}\fi}



\newtoks\auteurcourant      \auteurcourant={\hfil}
\newtoks\titrecourant       \titrecourant={\hfil}

\newtoks\hautpagetitre      \hautpagetitre={\hfil}
\newtoks\baspagetitre       \baspagetitre={\hfil}

\newtoks\hautpagegauche
\hautpagegauche={\eightpoint\rlap{\folio}\hfil\the\auteurcourant\hfil}
\newtoks\hautpagedroite
\hautpagedroite={\eightpoint\hfil\the\titrecourant\hfil\llap{\folio}}

\newtoks\baspagegauche      \baspagegauche={\hfil}
\newtoks\baspagedroite      \baspagedroite={\hfil}

\newif\ifpagetitre          \pagetitretrue

\def\raggedbottom{\topskip 10pt plus 36pt\r@ggedbottomtrue}

\def\hfl#1#2{\nospacedmath\smash{\mathop{\hbox to
12mm{\rightarrowfill}}\limits^{\scriptstyle#1}_{\scriptstyle#2}}}

\def\phfl#1#2{\nospacedmath\smash{\mathop{\hbox to
8mm{\rightarrowfill}}\limits^{\scriptstyle#1}_{\scriptstyle#2}}}

\def\vfl#1#2{\llap{$\scriptstyle#1$}\left\downarrow\vbox to
6mm{}\right.\rlap{$\scriptstyle#2$}}

\long\def\ctitre#1\endctitre{%
    \ifdim\lastskip<24pt\vskip-\lastskip\bigbreak\bigbreak\fi
  		\vbox{\parindent=0pt\leftskip=0pt plus 1fill
          \rightskip=\leftskip
          \parfillskip=0pt\bf#1\par}
    \bigskip\nobreak}

\mathcode`A="7041 \mathcode`B="7042 \mathcode`C="7043
\mathcode`D="7044
\mathcode`E="7045 \mathcode`F="7046 \mathcode`G="7047
\mathcode`H="7048
\mathcode`I="7049 \mathcode`J="704A \mathcode`K="704B
\mathcode`L="704C
\mathcode`M="704D \mathcode`N="704E \mathcode`O="704F
\mathcode`P="7050
\mathcode`Q="7051 \mathcode`R="7052 \mathcode`S="7053
\mathcode`T="7054
\mathcode`U="7055 \mathcode`V="7056 \mathcode`W="7057
\mathcode`X="7058
\mathcode`Y="7059 \mathcode`Z="705A

\def\spacedmath#1{\def\packedmath##1${\bgroup\mathsurround=0pt
##1\egroup$}%
\mathsurround#1
\everymath={\packedmath}\everydisplay={\mathsurround=0pt }}

\def\nospacedmath{\mathsurround=0pt \everymath={}\everydisplay={} }

\def\up#1{\raise 1ex\hbox{\smallf@nt#1}}


\def\cf{{\it cf}} 

\def\qed{\raise -2pt\hbox{\vrule\vbox to 10pt{\hrule width 4pt
                 \vfill\hrule}\vrule}}

\def\cqfd{\unskip\penalty 500\quad\vrule height  4pt depth 0pt width
4pt\medbreak}

\def\virg{\raise .4ex\hbox{,}}   


\def\date{\the\day\ \ifcase\month\or janvier\or f\'evrier\or mars\or
avril\or mai\or juin\or juillet\or ao\^ut\or septembre\or octobre\or
novembre\or d\'ecembre\fi \ {\old \the\year}}

\def\dateam{\ifcase\month\or January\or February\or March\or
April\or May\or June\or July\or August\or September\or October\or
November\or December\fi \ \the\day ,\ \the\year}

\def\crog{{\vrule height 2.57mm depth 0.85mm width 0.3mm}\kern -0.36mm
[}

\def\crod{]\kern -0.4mm{\vrule height 2.57mm depth 0.85mm
width 0.3 mm}}

\def\rond{\kern 1pt{\scriptstyle\circ}\kern 1pt}

\def\diagram#1{\def\normalbaselines{\baselineskip=0pt
\lineskip=10pt\lineskiplimit=1pt}   \matrix{#1}}

\def\av{abelian variety}
\def\avs {abelian varieties}

\def\pa{\S\kern.15em}

\def\P{\hbox{\bf P}}

\def\C{\hbox{\bf C}}

\def\cf{{\it cf.\/}}

\def\codim{\mathop{\rm codim}\nolimits}

\def\isom{\simeq}

\def\loc{{\it loc.cit.\/}}
\def\long{\mathop{\rm long}\nolimits}
\def\lra{\longrightarrow}
\def\llra{\nospacedmath\hbox to 10mm{\rightarrowfill}}
\def\lllra{\nospacedmath\hbox to 15mm{\rightarrowfill}}

\def\Pic{\mathop{\rm Pic}\nolimits}

\def\ra{\rightarrow}

\def\tx{\kern -1.5pt -}

\def\cc#1{\hfill\kern .7em#1\kern .7em\hfill}


\catcode`\@=12

\showboxbreadth=-1  \showboxdepth=-1



\baselineskip=13.8pt
\spacedmath{2pt}
\font\eightrm=cmr10 at 8pt
\overfullrule=0mm
\parskip=1.7mm
\def\ndeg{non-degenerate}
\def\gndeg{geometrically non-degenerate}
\def\iff{if and only if}
\def\div{\mathop{\rm div}\nolimits}

\null\smallskip
\ctitre  {\bf  FULTON-HANSEN AND BARTH-LEFSCHETZ THEOREMS FOR  SUBVARIETIES
OF ABELIAN VARIETIES}
\endctitre\smallskip
\smallskip
\centerline{Olivier {\pc DEBARRE}\footnote{(*)}{\rm Partially supported by
N.S.F. Grant DMS 94-00636 and the European HCM project ``Algebraic Geometry
in Europe" (AGE), Contract CHRXCT-940557.}}
\vskip 1cm
\ind Sommese showed that a large part of the geometry of a {\it smooth}
subvariety of a complex \av\ depends on ``how ample'' its normal bundle is
(see \S 1 for more details). Unfortunately, the only known way of measuring
this ampleness uses rather strong properties  of the ambiant \av .

\ind We show that a notion of non-degeneracy due to Ran is a good substitute
for ampleness of the normal bundle. It can be defined as follows: {\it an
irreducible subvariety $V$ of an \av\ $X$ is \gndeg\ if for any abelian
variety $Y$ quotient of $X$, the image of $V$ in $Y$ either is $Y$ or has
same dimension as $V$.} This property does not require $V$ to be smooth; for
smooth subvarieties, it is (strictly) weaker than ampleness of the normal
bundle.

\ind Our main result is a Fulton-Hansen type theorem for an irreducible
subvariety $V$ of an \av : the dimension of the ``secant variety'' of $V$
along a subvariety $S$ (defined as
$V-S$), and that of its ``tangential variety'' along $S$ (defined in the
smooth case as the union of the
 projectivized tangent spaces to $V$ at points of $S$, translated at the
origin) differ by $1$. Corollaries include a new proof of the finiteness of
the Gauss map and an estimate on the ampleness of the normal bundle of a
smooth \gndeg\ subvariety.

\ind We also complement Sommese's work with a new Barth-Lefschetz theorem for
subvarieties of
\avs\ whose proof is based on an idea of Schneider and Zintl. Let $C$ be a
smooth curve in an
\av\ $X$; we apply this result to give an estimate on the dimension of the
singular locus of
$C+\cdots +C$ in $ X$.

\ind We work over the field of complex numbers.
\bigskip

{\bf 1. Geometrically non-degenerate subvarieties}

\ind Recall ([S1]) that a line bundle $L$ on an irreducible projective
variety $V$ is $k\!$--ample if, for some $m>0$, the line bundle $L^m$ is
generated by its global sections and the fibers of the associated map
$\phi_{L^m}:V\ra
\P^N$ are all of dimension $\le k$. A vector bundle $E$ on $V$ is
$k\!$--ample if the line bundle ${\cal O}_{{\bf P}E^*}(1)$ is $k\!$--ample.
Ordinary ampleness coincide with
$0$--ampleness.

\ind Let $X$ be an \av\ and let $V$ be an irreducible variety with a morphism
$f:V\ra X$. Let
$V^0$ be the open set of smooth points of $V$ at which $f$ is unramified.
Define the normal bundle to $f$ as the vector bundle on $V^0$ quotient of
$f^*(TX)_{|V^0}$ by $TV^0$.

\ind For any $x\in X$, let $\tau_x$ be the translation by $x$. For any $v\in
V$, the differential of the map $\tau_{-f(v)}f$ at $v$ is a linear map
$T_vV\ra T_0X$, which we will simply denote by
$f_*$.
\bigskip {\pc PROPOSITION} 1.1.-- {\it Under the above assumptions, let $S$ be
a complete irreducible subvariety of $V^0$. The following properties are
equivalent

\ind{\rm (i)} the restriction to $S$ of the normal bundle to $f$ is
$k\!$--ample;

\ind{\rm (ii)} for any hyperplane $H$ in $T_0X$, the set $\{\ s\in S\bigm|
f_*(T_sV)\subset H\ \}$ has dimension $\le k$.}

{\bf Proof}. Let $N$ be the restriction to $S$ of the normal bundle to $f$
and let $\iota:\P N^*\ra\P f^*(T^*X)_{|S}$ be the canonical injection. The
morphism
$$\phi:\P N^* \buildrel {\iota}\over{\lra}\P f^*(T^*X)_{|S}\isom\P
T^*_0X\times S \buildrel {pr_1}\over{\lra}\P T^*_0X$$ satisfies
$$\phi^*{\cal O}_{{\bf P} T^*_0X}(1)=\iota^*{\cal O}_{{\bf P}
f^*(T^*X)_{|S}}(1)={\cal O}_{{\bf P}N^*}(1)\ .$$ \ind It follows that $N$ is
$k\!$--ample if and only if the fibers of $\phi$ have dimension $\le k$
([S1], prop. 1.7). The proposition follows, since the restriction of the
projection $\P N^* \ra S$ to any fiber of
$\phi$
 is injective.\cqfd

\ind When $X$ is simple, the normal bundle to any smooth subvariety of $X$ is
ample ([H]). More generally, the normal bundle to any smooth subvariety of
$X$ is $k\!$--ample, where $k$ is the maximum dimension of a proper abelian
subvariety of $X$  ([S1], prop. 1.20).

\ind Following Ran, we will say that a $d$\tx dimensional irreducible
subvariety $V$ of $X$ is {\it \gndeg } if the kernel of the restriction
$H^0(X,\Omega^d_X)\ra H^0(V_{\rm reg},\Omega^d_{V_{\rm reg}})$ contains no
non-zero decomposable forms. This property holds \iff\ for any abelian
variety $Y$ quotient of $X$, the image of $V$ in $Y$ either is $Y$ or has
same dimension as $V$ ([R1], lemma II.12).

{\it Examples.} 1) A divisor is \gndeg\ \iff\ it is ample; a curve is \gndeg\
\iff\ it generates $X$. Any \gndeg\ subvariety of positive dimension
generates $X$, but the converse is false in general. However, any irreducible
subvariety of a {\it simple} \av\ is \gndeg .

\ind 2) If $\ell$ is a polarization on $X$ and $V$ is an irreducible
subvariety of $X$ with class a rational multiple of $\ell^c$, it follows from
[R1], cor. II.2 and II.3 that $V$ is \ndeg\ in the sense of [R1], II, hence
\gndeg . In particular, the subvarieties $W_d(C)$ of the Jacobian of a curve
$C$ are \gndeg ; it can be checked  that their normal bundle is ample when
they are smooth (use prop. 1.1).
 \medskip
\ind We generalize this notion as follows. Let $k$ be a non-negative integer.
\medskip {\pc DEFINITION} 1.2.-- {\it An irreducible subvariety $V$ of an
\av\ $X$ is $k\!$--\gndeg\ \iff\ for any abelian variety $Y$ quotient of $X$,
the image of $V$ in $Y$ either is $Y$ or has dimension $\ge\dim (V)-k$.}
\bigskip {\pc PROPOSITION} 1.3.-- {\it In an \av , any smooth irreducible
subvariety with $k\!$--ample normal bundle is $k\!$--\gndeg .}

{\bf Proof}. Let $\pi:X\ra Y$ be a quotient of $X$ such that $\pi(V)\ne Y$.
The tangent spaces to $V$ along a smooth fiber of $\pi_{|V}$ are all
contained in a fixed hyperplane, hence general fibers of $\pi_{|V}$ have
dimension $\le k$ by prop. 1.1.\cqfd

\ind The converse is not true, as the construction sketched below shows, but
a partial converse will be obtained in 2.3. Roughly speaking, if $Y$ is a
quotient of $X$, and if the image $W$ of $V$ in $Y$ is not $Y$,
$k\!$--geometrical nondegeneracy requires that the general fibers of $V\ra W$
be of dimension $\le k$, whereas $k\!$--ampleness of the normal bundle
requires that every fiber of  $V\ra W$ be of dimension $\le k$.

\ind Let $L_E$ be an ample line bundle on an elliptic curve $E$, with
linearly independent sections $s_1,s_2$ defining a morphism $E\ra\P^1$ with
ramification points
$(e_1,1),\ldots,(e_4,1)\in\P^1$. Let $L_Y$ be an ample line bundle on a
simple abelian variety $Y$ of dimension $\ge 3$, with linearly independent
sections $t_1,t_2,t_3$ such that $\div (t_3)$,  $F=\div (t_1)\cap \div
(t_2)\cap \div (t_3)$ and $\div (e_it_1+t_2)\cap\div (t_3)$ are smooth  for
$i=1,\ldots,4$ (such a configuration can be constructed using results from
[D2]). Set
$X=E\times Y$ and define a subvariety $V$ of $X$ by the equations
$s_1t_1+s_2t_2=t_3=0$; then $V$ is smooth of codimension $2$, \gndeg , but
its normal bundle is not ample (for all $e\in E$ and $f\in F$, one has
$T_{(e,f)}V\subset T_f\bigl(\div (t_3)\bigr)$), only $1$--ample (cor. 2.3).

 \bigskip {\pc PROPOSITION} 1.4.-- {\it Let $X$ be an \av\ and let $V$ and
$W$ be irreducible subvarieties of
$X$. Define a morphism $\phi:V^r\times W\ra X^r$ by $\phi(v_1,
\ldots,v_r,w)=(v_1-w,\ldots,v_r-w)$. If $V$ is $k\!$--\gndeg ,
$$\dim \phi(V^r\times W)\ge \min \bigl( r\dim(X), r\dim(V)+\dim(W)-k\bigr)\
.$$} {\bf Proof}. Assume first $r\dim(V)+\dim(W)-k\ge r\dim(X)$. Let
$\pi:X\ra X/K$ be a quotient of $X$. I claim that
$r\dim\pi(V)+\dim \pi(W)\ge r\dim (X/K)$. If $\pi(V)=X/K$, this is obvious;
otherwise, we have
$\dim \pi(V)\ge \dim (V)-k_0$, where $k_0=\min \bigl( k,\dim(K)\bigr)$, hence
$$\eqalign{r\dim \pi(V)+\dim \pi(W)& \ge r(\dim (V)-k_0)+\dim (W)-\dim (K)\cr
& \ge r\dim(X)+k-rk_0-\dim(K)\ge r\dim (X/K)\ .\cr}$$
\ind It follows that $(V,\ldots,V,W)$ (where $V$ is repeated $r$ times) fills
up $X$ in the sense of [D1], (1.10);  th. 2.1 of \loc\ then implies that
$\phi$ is onto.

\ind Assume now $s=r\codim(V) -\dim(W)+k>0$; let $C$ be an irreducible curve
in $X$ that
 generates $X$. Let $W'$ be the sum of $W$ and $s$ copies of $C$; then
$r\dim(V)+\dim(W')-k= r\dim(X)$ and the first case shows that the sum of the
image of $\phi$ and $s$  curves is $X^r$. The proposition follows.\cqfd
\medskip
\ind We obtain a nice characterization of  $k\!$--\gndeg\ varieties.
\medskip {\pc COROLLARY} 1.5.-- {\it An irreducible subvariety $V$ of an
abelian variety $X$ is $k\!$--\gndeg\ \iff\ it meets any subvariety of $X$ of
dimension
$\ge\codim(V)+k$.}

{\bf Proof}. Assume that $V$ meets any subvariety of $X$ of dimension
$\ge\codim(V)+k$ and let $\pi:X\ra Y$ be a quotient of $X$. If
$\pi(V)\ne Y$, there exists a subvariety $W$ of $Y$ of dimension
$\dim(Y)-\dim\pi(V)-1$ that does not meet $\pi(V)$. Since $V$ does not meet
$\pi^{-1}(W)$,
$$\codim (V)+k>\dim \pi^{-1}(W)=\dim(X)-\dim\pi(V)-1$$ hence $\dim \pi(V)\ge
\dim(V)-k$ and $V$ is $k\!$--\gndeg . Conversely, assume $V$ is
$k\!$--\gndeg; let $W$ be an irreducible subvariety of $X$ of dimension
$\ge\codim(V)+k$. Proposition 1.4 shows that $V-W=X$, hence $V$ meets
$W$.\cqfd

 \bigskip

{\bf 2. A Fulton-Hansen-type result}

\ind Fulton and Hansen proved in [FH] (\cf\ also [FL1], [Z1], [Z2]) a
beautiful result that relates the dimension of the tangent variety and that
of the secant variety of a subvariety of a projective space. We prove an
analogous result for a subvariety of an \av .

\ind Let $X$ be an \av\ and let $V$ be a variety with a morphism $f:V\ra X$.
Recall that $f$ is unramified along a subvariety $S$ of $V$ if $\Delta_S=\{
(v,s)\in V\times S\mid v=s \}$ is an open subscheme of $V\times_XS$.
Following [FL1], we will say that $f$ is {\it weakly unramified} along  $S$
if $\Delta_S$ is a connected component of $V\times_XS$, ignoring scheme
structures. In that case, if $p:V\times S\ra X$ is the morphism defined by
$p(v,s)=f(v)-f(s)$ and $\epsilon:\tilde X\ra X$ is the blow-up of the origin,
there exists a commutative diagram
$$\diagram{\tilde Y& \hfl{\tilde p}{}& \tilde X\cr
\vfl{\alpha}{}&&\vfl{\epsilon}{}\cr V\times S& \hfl{p}{}& X\cr}$$ where
$\alpha$ is the blow-up of $V\times_XS$. Let $E$ be the exceptional divisor
above
$\Delta_S\i V\times S$ and set $T(V,S)=\tilde p(E)$. It is a subscheme of $\P
T_0X$ contained in
$\bigcup_{s\in S}\P T_sV$ and equal to the latter when $V$ is smooth along
$S$ and $f$ is unramified along $S$. Loosely speaking, $T(V,S)$ is the set of
limits in $\tilde X$ of
$\bigl( f(v)-f(s)\bigr)$, as $v\in V$ and $s\in S$ converge to the same
point. Obviously, $\dim T(V,S)<\dim\bigl( f(V)-f(S)\bigr)$.   \bigskip {\pc
THEOREM} 2.1.-- {\it Let $X$ be an \av\ and let $V$ be an irreducible
projective variety with a morphism $f:V\ra X$. Let $S$ be a complete
irreducible subvariety of $V$ along which $f$ is weakly unramified. Then
$\dim T(V,S)=\dim\bigl( f(V)-f(S)\bigr) -1$.}

\ind We begin with a lemma.
\smallskip {\pc LEMMA} 2.2.-- {\it Let $C$ be an irreducible projective curve
with a morphism $g:C\ra X$ such that $g(C)$ is a smooth curve  through the
origin. Assume that $g$ is unramified at some point $c_0\in C$ with
$g(c_0)=0$ and that $\P T_0g(C)\notin T(V,S)$. The morphism $h:V\times C\ra
X$ defined by $h(v,c)=f(v)-g(c)$ is weakly unramified along $S\times\{ c_0\}$
and $T(V\times C,S\times\{ c_0\} )$ is contained in the cone over
$T(V,S)$ with vertex $\P T_0g(C)$.}

\ind One can prove that $T(V\times C,S\times\{ c_0\} )$ is actually equal to
the cone.

{\bf Proof}. Let $\Gamma$ be a smooth irreducible curve, let $\gamma_0$ be a
point on $\Gamma$ and let
$q=(q_1,q_2,q_3):\Gamma\ra (V\times C)\times_XS$ be a morphism with
$q(\gamma_0)=(s_0,c_0,s_0)$. We need to prove that
$q(\Gamma)\subset\Delta'_S$, where $\Delta'_S=\{(s,c_0,s)\mid s\in S\}$; since
$\Delta_S$ is a connected component of $V\times_XS$, it suffices to show that
$q_2$ is constant. Suppose the contrary; then $(q_1,q_3)$ lifts to a morphism
$\tilde q_{13}:\Gamma\ra
\tilde Y$ and $g$ to a morphism $\tilde g:C\ra\tilde X$. Since
$p(q_1,q_3)=gq_2$, one has
$\tilde p\tilde q_{13}=\tilde gq_2$ hence $\tilde g(c_0)=\tilde p\bigl(\tilde
q_{13}(\gamma_0)\bigr)\in T(V,S)$. This contradicts the hypothesis since
$\tilde g(c_0)$ is the point
$\P T_0g(C)$ of $\P T_0X$. This proves the first part of the lemma.

\ind The second part is similar: let $\tilde Z\ra (V\times C)\times S$ be the
blow-up of $(V\times C)\times_XS$, let $\Gamma$ be a smooth irreducible curve
with a point $\gamma_0\in \Gamma$ and let
$\tilde q:\Gamma\ra \tilde Z$ be a morphism such that $\tilde q(\gamma_0)$ is
in the exceptional divisor above $\Delta'_S$. Write $q=\alpha\tilde
q=(q_1,q_2,q_3)$ and keep the same notation as above. Then $pq(\Gamma)$ is
contained in the surface $pq_{13}(\Gamma)-g(C)$ hence
$\tilde p\tilde q(\gamma_0)$ belongs to the line in $\P T_0X$ through $\tilde
p\bigl(\tilde q_{13}(\gamma_0)\bigr)$ and  $\tilde g(c_0)=\P T_0g(C)$. This
proves the lemma.\cqfd

{\bf Proof of the theorem}. We proceed by induction on the codimension of
$f(V)-f(S)$. Assume $f(V)-f(S)=X$; if $T(V,S)\ne\P T_0X$, pick a point
$u\notin T(V,S)$ and a smooth projective curve $C'$ in $X$ tangent to $u$ at
$0$, and such that the restriction induces an injection
$\Pic^0(X)\ra\Pic^0(C')$. Let $C$ be a smooth curve with a connected ramified
double cover $g:C\ra C'$ unramified at a point $c_0$ above $0$; the map
$\Pic^0(C')\ra\Pic^0(C) $ induced by $g$ is injective.

\ind Since $p$ is surjective, $C'$ generates $X$ and $C$ is smooth, th. 3.6
of [D1] implies that $(V\times S)\times_XC$ is connected. If $h:V\times C\ra
X$ is defined by
$h(v,c)=f(v)-g(c)$, it follows that $(V\times C)\times_X S$ is also
connected.  On the other hand, the lemma implies that the set $\{
(s,c_0,s)\mid s\in S \}$ is a connected component of, hence is equal to,
$(V\times C)\times_X S$. It follows that $h^{-1}\bigl( f(S)\bigr) =S\times\{
c_0\}$. Since $g^{-1}(0)$ consists of $2$ distinct points, this is absurd,
hence $T(V,S)=\P T_0X$ and the theorem holds in this case.

\ind Assume now $f(V)-f(S)\ne X$. Take a curve $C'$ as above; by the lemma,
the morphism $f':V\times C'\ra X$ defined by $f'(v,c')=f(v)+c'$ is weakly
unramified along $S\times\{ 0\}$, and
$\dim T(V\times C',S\times\{ 0\} )\le\dim T(V,S)+1$. It follows from the
induction hypothesis that
$$\dim T(V,S)\ge\dim \bigl( f(V)+C'-f(S)\bigr) -2=\dim \bigl( f(V)-f(S)\bigr)
-1\ ,$$ which proves the theorem.\cqfd
\ind The following corollary provides a partial converse to prop. 1.3.
\bigskip {\pc COROLLARY} 2.3.-- {\it Let $X$ be an \av\ of dimension $n$ and
let $V$ be an irreducible projective variety of dimension $d$ with a morphism
$f:V\ra X$ such that $f(V)$ is $k\!$--\gndeg . Let $V^0$ be the open set of
smooth points of $V$ at which $f$ is unramified. The restriction of the
normal bundle to $f$ to any complete irreducible subvariety $S$ of $V^0$ is
$(n-d-1+k)\!$--ample.}

{\bf Proof}. By prop. 1.1, we must show that for any hyperplane $H$ in
$T_0X$, any irreducible component $S_H$ of $\{ s\in S\mid f_*(T_vV)\subset
H\}$ has dimension $\le n-d-1+k$. But $T(V,S_H)$ is contained in $H$ and the
theorem gives $f(V)-f(S_H)\ne X$. Since $f$ is unramified along $S_H$ and
$f(V)$ is $k\!$--\gndeg , prop. 1.4 implies that
$f(V)-f(S_H)$ has dimension $\ge d+\dim(S_H)-k$; this proves the
corollary.\cqfd

\ind It should be noted that the corollary also follows from the main result
of [Z3] (cor.~1),
 whose proof is unfortunately so sketchy (to say the least) that I could not
understand it.
 \bigskip
 {\pc COROLLARY} 2.4.-- {\it Let $X$ be an \av\ and let $V$ be an irreducible
projective variety with a morphism $f:V\ra X$. Let $L$ be a linear subspace
of $T_0X$ and let $S$ be a complete irreducible subvariety of $V$ along which
$f$ is unramified. Assume that $\dim\bigl( f_*(T_sV)\cap L\bigr) <m$ for all
$s\in S$,
 and let $\Delta_{f(S)}$ be the small diagonal in $f(S)^m$. Then
 $$\dim \bigl( f(V)^m-\Delta_{f(S)}\bigr)< m\dim(X)-\dim(L)+m\ .$$
\ind In particular, if $m\le \dim(L)$ and $f(V)$ is $k\!$--\gndeg ,
$$m\dim  (V)+\dim (S)< m\dim(X)-\dim(L)+m+k\ .$$} {\bf Proof}. Let
$r=\dim(L)$; the variety $N=\{\ [t_1,\ldots,t_m]\in \P(L^m)\mid
t_1\wedge\cdots\wedge t_m=0\}$ has codimension $r-m+1$ in $\P(L^m)$. Consider
the morphism $f^m:V^m\ra X^m$ and the subvariety $\Delta_S$ of $V^m$. The
hypothesis imply that in $\P (T_0X^m)$, the intersection of $T(V^m,\Delta_S)$
and $\P (L^m)$ is contained in  $N$. It follows that $$\eqalign{\dim
T(V^m,\Delta_S)&\le \dim (N) +\dim \P (T_0X^m)-\dim\P (L^m)\cr &= \dim \P
(T_0X^m)-(r-m+1)\ .\cr}$$
\ind The first inequality of the corollary follows from th. 2.1, and the
second from prop. 1.4.\cqfd \bigskip  {\bf 3. Applications to the Gauss map}

\ind We keep the same setting:  $X$ is an \av\ and $V$ an irreducible
projective variety of dimension $d$ with a morphism $f:V\ra X$. Let $V^0$ be
the open set of smooth points of $V$ at which $f$ is unramified; define the
{\it Gauss map} $\gamma:V^0\ra G(d,T_0X)$ by
$\gamma(v)=f_*(T_vV)$. The following result was first proved by Ran ([R2]),
and by Abramovich ([A]) in all characteristics.

\medskip {\pc PROPOSITION} 3.1.-- {\it Let $X$ be an \av\ and let $V$ be an
irreducible projective variety with a morphism $f:V\ra X$. If $S$ is a
complete irreducible variety contained in a fiber of the Gauss map,
$f(V)$ is stable by translation by the abelian variety generated by $f(S)$. In
particular, the Gauss map of a smooth projective subvariety of $X$ invariant
by translation by no non-zero abelian subvariety of $X$  is finite.}

{\bf Proof}. Under the hypothesis of the proposition, $T(V,S)$ has dimension
$\dim(V)-1$;  th. 2.1 implies $f(V)-f(S)=f(V)$, hence the proposition.\cqfd
\medskip
\ind For any linear subspace $L$ of $T_0X$  and any integer $m\le
\dim(L)$, let $\Sigma_{L,m}$ be the Schubert variety $\{ M\in G(d,T_0X)\mid
\dim (L\cap M)\ge m \}$; its codimension in $G(d,T_0X)$ is
$m(\codim(L)-d+m)$.
\medskip   {\pc PROPOSITION} 3.2.-- {\it Let $X$ be an
\av\ and let $V$ be an irreducible projective variety of dimension $d$ with a
morphism $f:V\ra X$ such that $f(V)$ is $k\!$--\gndeg . Let $\gamma:V^0\ra
G(d,T_0X)$ be the Gauss map, let
$L$ be a linear subspace of $T_0X$ and let $m$ be an integer $\le\dim(L)$.
Any complete subvariety $S$ of $V^0$ of dimension $\ge \codim
\Sigma_{L,m}+(m-1)\bigl(\dim(L)-m\bigr) +k$ meets
$\gamma^{-1}(\Sigma_{L,m})$.}

{\bf Proof}. Apply cor. 2.4.\cqfd

\ind The hypothesis could probably be weakened to
$\dim(S)\ge \codim
\Sigma_{L,m}+k$ (see next proposition); the proposition gives that for
$m=1$ or $\dim(L)$. The corresponding Schubert varieties are $\Sigma_{L,1}=\{
M\in G(d,T_0X)\mid L\cap M\ne 0 \}$ and $\Sigma_{L,\dim(L)}=\{ M\in
G(d,T_0X)\mid  L\subset M \}$.

\ind More generally, a result of Fulton and Lazarsfeld imposes strong
restrictions on the image of the Gauss map of smooth subvarieties with ample
normal bundle which I believe should also hold for \gndeg\ subvarieties.
\medskip

{\pc PROPOSITION} 3.3.-- {\it Let $X$ be an \av\ and let $V$ be a smooth
irreducible projective variety of dimension $d$ with an unramified morphism
$f:V\ra X$ and Gauss map
$\gamma:V\ra G(d,T_0X)$. Assume that the normal bundle to $f$ is ample; any
subvariety $S$ of $\gamma(V)$ meets any subvariety of $G(d,T_0X)$ of
codimension $\le \dim (S)$.}

{\bf Proof}. If $Q$ is the universal quotient bundle on $G(d,T_0X)$, the
pull-back $\gamma^*(Q)$ is isomorphic to the normal bundle $f^*TX/TV$, hence
is ample. It follows from [FL2] that for each Schubert variety
$\Sigma_{\lambda}$ of codimension $m$ in $G(d,T_0X)$ and each irreducible
subvariety $S$ of $V$ of dimension $m$, one has
$\int_S\gamma^*[\Sigma_{\lambda}]>0$. Now the class of any irreducible
subvariety $Z$ of
$G(d,T_0X)$ of codimension $m$ is a linear combination with non-negative
coefficients (not all zero) of the  Schubert classes; this implies
$\int_S\gamma^*[Z]>0$, hence $S\cap
\gamma^{-1}(Z)\ne \emptyset$.\cqfd

\medskip
\ind Regarding the Gauss map of a smooth subvariety of an \av ,  Sommese and
Van de Ven also proved in [SV] a  strong result for higher relative homotopy
groups of pull-backs of {\it smooth} subvarieties of the Grassmannian.

\bigskip {\bf 4. A Barth-Lefschetz-type result}

\ind Sommese has obtained very complete results on the homotopy groups of
{\it smooth} subvarieties of an \av . For example, he proved in [S2] that if
$V$ is a smooth subvariety of dimension $d$ of an \av\ $X$,  with
$k\!$--ample normal bundle, $\pi_q(X,V)=0$ for $q\le 2d-n-k+1$. For arbitrary
subvarieties, we have the  following:\medskip {\pc THEOREM} 4.1.-- {\it Let
$X$ be an \av\  and let $V$ be a $k\!$--\gndeg\ normal subvariety of $X$ of
dimension $>{1\over 2}\bigl(\dim(X)+k\bigr)$. Then $\pi_1^{\rm alg}(V)\isom
\pi_1^{\rm alg}(X)$.}

{\bf Proof}. The case $k=0$ is cor. 4.2 of [D1]. The general case is similar,
since the hypothesis implies that the pair $(V,V)$ satisfies condition $(*)$
of [D1].\cqfd

\ind Going back to smooth subvarieties, I will give an elementary proof of (a
slight improvement of) the cohomological version of Sommese's theorem, based
on the following vanishing theorem ([LP]) and the ideas of [SZ].
\medskip {\pc VANISHING} {\pc THEOREM} 4.2 (Le Potier, Sommese).-- {\it Let
$E$ be a $k\!$--ample rank $r$ vector bundle on a smooth irreducible
projective variety $V$ of dimension $d$. Then
$$H^q(V,E^*\otimes\Omega^p_V)=0\qquad{\it for}\qquad p+q\le d-r-k\ .$$}
\ind Recall also the following elementary lemma from [SZ]: \medskip {\pc
LEMMA} 4.3.-- {\it Let $0\ra F\ra E_0\ra E_1\ra\ldots\ra E_k\ra 0$ be an exact
sequence of sheaves on a scheme $V$. Assume $H^s(V,E_i)=0$ for
$0\le i<k$ and $s\le q$; then $H^q(V,F)\isom H^{q-k}(V,E_k)$.}
\medskip {\pc THEOREM} 4.4.-- {\it Let
$V$ be a smooth irreducible subvariety of dimension $d$ of an abelian
$n\!$--fold $X$ and let ${\cal L}$ be a nef line bundle on $V$. Assume that
the normal bundle $N$ of $V$ in $X$ is a direct sum $\oplus N_i$, where
$N_i$ is $k_i\!$--ample of rank $r_i$. For $j>0$}, $$H^q(V,S^j
N^*\otimes{\cal L}^{-1})=0\qquad{\it for}\qquad q\le d-\max (r_i+k_i)\ .$$
{\bf Proof}. Since $S^jN^*$ is a direct summand of $S^{j-1}N^*\otimes N^*$,
it is enough to show, by induction on $j$, that $H^q(V,S^j N^*\otimes
N_i^*\otimes{\cal L}^{-1})$ vanishes for $j\ge 0$ and $q\le d-r_i-k_i$. Since
$N_i\otimes{\cal L}$ is $k_i\!$--ample, the case $j=0$ follows from Le
Potier's theorem. For $j\ge 1$, tensor the exact sequence
$$0\ra S^jN^*\ra
S^{j-1}N^*\otimes\Omega^1_{X|V}\ra\ldots\ra\Omega^j_{X|V}\ra\Omega^j_V\ra 0$$
by $N_i^*\otimes{\cal L}^{-1}$. Since $\Omega_X^1$ is trivial, the induction
hypothesis and the lemma give $H^q(V,S^jN^*\otimes N_i^*\otimes{\cal
L}^{-1})\isom H^{q-j}(V,\Omega^j_V\otimes N_i^*\otimes{\cal L}^{-1})$, and
this group vanishes for\break $q\le d-r_i-k_i$ by Le Potier's theorem.\cqfd

\ind We are now ready to prove our version of Sommese's result.
\bigskip {\pc THEOREM} 4.5.-- {\it Let
$V$ be a smooth irreducible subvariety of dimension $d$ of an abelian
$n\!$--fold $X$. Assume that its normal bundle a direct sum $\oplus N_i$,
where $N_i$ is
$k_i\!$--ample of rank $r_i$. Then

\ind {\rm a)} $H^q(X,V;\C )=0$ for $q\le d-\max (r_i+k_i)+1${\rm ;

\ind  b)} for all nonzero elements $P$ of $\Pic ^0(V)$, the cohomology groups
$H^q(V,P)$ vanish for $q\le d-\max (r_i+k_i)$.}
\medskip {\it Remarks} 4.6. 1) It is likely that a) should hold for
cohomology with integral coefficients.

\ind 2) If the normal bundle is $k\!$--ample, we get $H^q(X,V;\C )=0$ for
$q\le 2d-n-k+1$. If the normal bundle is a sum of ample line bundles,
$H^q(X,V;\C )=0$ for $q\le d$; in particular, the restriction
$H^0(X,\Omega_X^d)\ra H^0(V,\Omega_V^d)$ is injective and $V$ is
non-degenerate in the sense of [R1], II, hence also \gndeg .

\ind 3) By [GL], $H^q(V,P)=0$  for $P$ outside of a subset of codimension
$\ge d-q$ of
$\Pic ^0(V)$. By [S], this subset is a union of translates of abelian
subvarieties of $X$ by torsion points.

\medskip {\bf Proof of the theorem}. For a), it is enough by Hodge theory to
study the maps
$$H^i(X,\Omega^j_X)\lra H^i(V,\Omega^j_{X|V})\buildrel\psi\over{\lra}
H^i(V,\Omega^j_V)\ .$$\ind Since $\Omega_X^j$ is trivial, we only need look
at $\phi:H^i(X,{\cal O}_X)\ra H^i(V,{\cal O}_V)$ and $\psi$. We begin with
$\psi$. We may assume $j>0$. Let $M_j$ be the kernel of the surjection
$\Omega^j_{X|V}\ra \Omega^j_V\ra 0$. The long exact sequence of th. 4.4 gives
$$0\ra S^jN^*\ra S^{j-1}N^*\otimes\Omega^1_{X|V}\ra\ldots\ra
N^*\otimes\Omega^{j-1}_{X|V}\ra M_j\ra 0\ .$$
\ind The lemma and the theorem then yield
 $$H^i(V,M_j)\isom H^{i+j-1}(V,S^jN^*)=0$$ for
$i+j-1\le d-\max (k_i+r_i)$, since $j>0$. This implies that $\psi$ has the
required properties.

\ind For $i=0$, the map $\psi$ is $H^0(X,\Omega_X^j)\ra H^0(V,\Omega_V^j)$.
By Hodge symmetry, this proves that $\phi$ also has the required properties,
hence the first point.

\ind For b), we may assume $d-\max (k_i+r_i)\ge 1$, in which case the first
point implies $\Pic^0(X)\isom \Pic^0(V)$. Let $P\in\Pic^0(X)$ be nonzero; the
same proof as above yields $H^0(V,\Omega_V^q\otimes P_{|V})=0$ for $q\le
d-\max (k_i+r_i)$. The theorem follows from the existence of an anti-linear
isomorphism $H^0(V,\Omega_V^q\otimes P_{|V})\isom H^q(V, P_{|V}^*)$
([GL]).\cqfd

\ind I will end this section with an amusing consequence of th. 4.5. If $C$
is a curve in an \av\
$X$, write $C_d$ for the subvariety $C+\cdots +C$ ($d$ times) of $X$. Recall
that if $C$ is general of genus $n$ and $d<n$, the singular locus of
$C_d=W_d(C)$ in the Jacobian $JC$ has dimension $2d-n-2$.

 {\pc PROPOSITION} 4.7.-- {\it Let $X$ be an abelian variety of dimension $n$
and let $C$ be a smooth irreducible curve in $X$. Assume that $C$  generates
$X$ and that its Gauss map is birational onto its image. Then, for $d< n$,
the singular locus of $C_d$ has dimension $\ge 2d-n-1$ unless $X$ is
isomorphic to the Jacobian of $C$ and $C$ is canonically embedded in $X$.}

{\bf Proof}. Let $\gamma:C_{\rm reg}\ra \P T_0X$ be the Gauss map and let
$\pi:C^{(d)}\ra C_d$ be the sum map. The image of the differential of $\pi$
at the point
$(c_1{\scriptscriptstyle\bullet}\ldots{\scriptscriptstyle\bullet}c_d)$ is the
linear subspace of $T_0X$ generated by $\gamma (c_1),\ldots,\gamma (c_d)$.
Since $C$ generates $X$, the curve $\gamma(C)$ is non-degenerate; it follows
that for
$c_1,\ldots,c_d$ general, the points $\gamma(c_1),\ldots,\gamma(c_d)$ span a
$(d-1)$--plane whose intersection with the curve $\gamma(C)$ consists only of
these points. Thus $\pi$ is birational. Moreover, if $x=c_1+\cdots +c_d$ is
smooth on
$C_d$, then $\gamma(c_i)\in \tau^*_x(T_xC_d)\cap\gamma(C)$ hence
$\pi^{-1}(x)$ is finite. By Zariski's Main Theorem, $\pi$ induces an
isomorphism between $\pi^{-1}\bigl( (C_d)_{\rm reg}\bigr)$ and $(C_d)_{\rm
reg}$.

\ind Let $s$ be the dimension of the singular locus of $C_d$ and assume
$-1\le s\le 2d-n-2$. Let $L$ be a very ample line bundle on $X$; the
intersection $W$ of $C_d$ with $(s+1)$ general elements of $|L|$ is smooth of
dimension $\ge 2$ and contained in
$(C_d)_{\rm reg}$. If $H$ is a hyperplane in $T_0X$ and $x=c_1+\cdots +c_d\in
W$, the inclusion
$T_xC_d\subset H$ implies  $\gamma(c_i)\in \P H\cap\gamma(C)$;  the
restriction of $N_{C_d/X}$ to $W$ is ample by prop. 1.1. Since $N_{W/X}$ is
the direct sum of this restriction and of $(s+1)$ copies of $L$,  the
restriction $H^1(X,{\cal O}_X)\ra H^1(W,{\cal O}_W)$ is bijective by th. 4.5.
On the other hand, the line bundle $\pi^*L$ is nef and big on
$C^{(d)}$, hence the Kawamata-Viehweg vanishing theorem ([K], [V]) implies
$$H^1(C^{(d)},{\cal O}_{C^{(d)}})\subset H^1(\pi^{-1}(W),{\cal
O}_{\pi^{-1}(W)})\isom H^1(W,{\cal O}_W)\ .$$
\ind Since $H^1(C,{\cal O}_C)\isom H^1(C^{(d)},{\cal O}_{C^{(d)}})$ ([M]), we
get
$h^1(C,{\cal O}_C)\le h^1(X,{\cal O}_X)$ and there must be equality because
$C$ generates $X$. Thus, the inclusion $C\subset X$ factors through an
isogeny $\phi:JC\ra X$. Since $\pi$ is birational, the inverse image
$\phi^{-1}(C_d)$ is the union of
$\deg (\phi)$ translates of $W_d(C)$. But any two translates of $W_d(C)$ meet
along a locus of dimension $\ge 2d-n> s$, hence $\phi$ is an isomorphism.\cqfd

\bigskip
\centerline{\pc REFERENCES}
\medskip\baselineskip=11.5pt

\hangindent=1cm [A] Abramovich, D., {\it Subvarieties of semiabelian
varieties},  Comp. Math. {\bf 90} (1994), 37--52.

\hangindent=1cm [D1] Debarre, O., {\it Th\'eor\`emes de connexit\'e et
vari\'et\'es ab\'eliennes}, Am. J. of Math. {\bf 117} (1995), 1--19.

\hangindent=1cm [D2] Debarre, O., {\it Sur les vari\'et\'es ab\'eliennes dont
le diviseur th\^eta est singulier en codimension $3$}, Duke Math. J. {\bf 56}
(1988), 221--273.

\hangindent=1cm [FH] Fulton, W., Hansen, J., {\it A connectedness theorem for
projective varieties, with applications to intersections and singularities of
mappings}, Ann. of Math. {\bf 110} (1979), 159--166.

\hangindent=1cm [FL1] Fulton, W., Lazarsfeld, R., Connectivity and its
Applications in Algebraic Geometry, in {\it Algebraic Geometry},  Proceedings
of the Midwest Algebraic Geometry Conference, Chicago 1980, Springer Lecture
Notes {\bf 862}.

\hangindent=1cm [FL2] Fulton, W., Lazarsfeld, R., {\it Positive polynomials
for ample vector bundles}, Ann. of Math. {\bf 118} (1983), 35--60.

\hangindent=1cm [GL] Green, M., Lazarsfeld, R., {\it Deformation Theory,
Generic Vanishing Theorems and Some Conjectures of Enriques, Catanese and
Beauville},  Invent. Math. {\bf 90} (1987), 389--407.

\hangindent=1cm [H] Hartshorne, R., {\it Ample Vector Bundles on Curves},
Nagoya Math. J. {\bf 43} (1971), 73--89.

\hangindent=1cm [K] Kawamata, Y., {\it A Generalisation of
Kodaira-Ramanujan's Vanishing Theorem},  Math. Ann. {\bf 261} (1982), 43--46.

\hangindent=1cm [LP] Le Potier, J., {\it Annulation de la cohomologie \`a
valeurs dans un fibr\'e vectoriel holomorphe positif de rang quelconque},
Math. Ann. {\bf 218} (1975), 35--53.

\hangindent=1cm [M] Macdonald, I., {\it Symmetric Products of an Algebraic
Curve},  Topology {\bf 1} (1962), 319--343.

\hangindent=1cm [R1] Ran, Z., {\it On subvarieties of abelian varieties},
Invent. Math. {\bf  62} (1981), 459--479.

\hangindent=1cm [R2] Ran, Z., {\it The structure of Gauss-like maps},  Comp.
Math. {\bf 52} (1984), 171--177.

\hangindent=1cm [SZ] Schneider, M., Zintl, J., {\it The theorem of
Barth-Lefschetz as a consequence of Le Potier's vanishing theorem},  Manuscr.
Math. {\bf 80} (1993), 259--263.

\hangindent=1cm  [S] Simpson, C., {\it Subspaces of Moduli Spaces of Rank One
Local Systems}, Ann. Sc. Ecole Norm. Sup. {\bf 26} (1993), 361--401.

\hangindent=1cm  [S1] Sommese, A., {\it Submanifolds of Abelian Varieties},
Math. Ann. {\bf 233} (1978), 229--256.

\hangindent=1cm  [S2] Sommese, A., {\it Complex Subspaces of Homogeneous
Complex Manifolds II. Homotopy Results},  Nagoya Math. J. {\bf 86} (1982),
101--129.

\hangindent=1cm  [SV] Sommese, A., Van de Ven, A., {\it Homotopy groups of
pullbacks of varieties}, Nagoya Math. J. {\bf 102} (1986), 79--90.

\hangindent=1cm  [V] Viehweg, E., {\it Vanishing Theorems},  J. reine Angew.
Math. {\bf 335} (1982), 1--8.

\hangindent=1cm  [Z1] Zak, F., {\it Structure of Gauss Maps}, Funct. Anal.
Appl. {\bf 21} (1987), 32--41.

\hangindent=1cm  [Z2] Zak, F., {\it Tangents and Secants of Algebraic
Varieties}, Transl. Math. Mono.  {\bf 127}, A.M.S., Providence, 1993.

\hangindent=1cm  [Z3] Zak, F., {\it Theorems on Tangents to Subvarieties of
Complex Tori}, Math. Notes Acad. Sc. USSR {\bf 43} (1988), 219--225.
\vskip 6mm
\hbox to 72mm{\hrulefill}\parskip=0cm {\pc OLIVIER} {\pc DEBARRE}

 IRMA -- {\pc MATH\'EMATIQUE} --  CNRS URA 001

{\pc UNIVERSIT\'E} {\pc LOUIS} {\pc PASTEUR}

7, rue Ren\'e Descartes

67084 {\pc STRASBOURG} {\pc CEDEX} -- {\pc FRANCE}

e-mail: debarre@math.u-strasbg.fr

\vfill\supereject\null\vskip 1cm

\ctitre  {\bf  FULTON-HANSEN AND BARTH-LEFSCHETZ THEOREMS FOR SUBVARIETIES OF
ABELIAN VARIETIES}
\endctitre\smallskip
\smallskip
\centerline{Olivier {\pc DEBARRE}}
\vskip 2cm {\bf Abstract}:
 We prove the following Fulton-Hansen type result for an irreducible
subvariety $V$ of an \av\
$X$: the dimension of the ``secant variety'' of $V$ along a subvariety $S$
(defined as
$V-S$), and that of its ``tangential variety'' along $S$ (defined in the
smooth case as the union of the
 projectivized tangent spaces to $V$ at points of $S$, translated at the
origin) differ by $1$. Corollaries include a new proof of the finiteness of
the Gauss map and an estimate on the ampleness of the normal bundle, for
certain smooth subvarieties of $X$. We also prove, using ideas of Schneider
and Zintl, a new Barth-Lefschetz theorem for smooth subvarieties of $X$. Let
$C$ be a smooth curve in $X$; we apply this result to give an estimate on the
dimension of the singular locus of
$C+\cdots +C$ in $X$.

\bye